\documentclass{amsart}
\usepackage{natbib}

\def\cal{\mathcal}
\newtheorem{theorem}{Theorem}[section]

\theoremstyle{definition}

\theoremstyle{remark}
\newtheorem{remark}[theorem]{Remark}

\numberwithin{equation}{section} \font\teneurm=eurm10
\font\seveneurm=eurm7 \font\fiveeurm=eurm5
\newfam\eurmfam
\textfont\eurmfam=\teneurm \scriptfont\eurmfam=\seveneurm
\scriptscriptfont\eurmfam=\fiveeurm

 \font\teneusm=eusm10 \font\seveneusm=eusm7 \font\fiveeusm=eusm5
\newfam\eusmfam
\textfont\eusmfam=\teneusm \scriptfont\eusmfam=\seveneusm
\scriptscriptfont\eusmfam=\fiveeusm
\def\eusm#1{{\fam\eusmfam\relax#1}}

\font\tencmmib=cmmib10 \skewchar\tencmmib='177
\font\sevencmmib=cmmib7 \skewchar\sevencmmib='177
\font\fivecmmib=cmmib5 \skewchar\fivecmmib='177
\newfam\cmmibfam
\textfont\cmmibfam=\tencmmib \scriptfont\cmmibfam=\sevencmmib
\scriptscriptfont\cmmibfam=\fivecmmib
\def\cmmib#1{{\fam\cmmibfam\relax#1}}

\def\Z{{\Bbb{Z}}}

\def\C{\Bbb{C}}
\def\M{{\cal M}}

\begin{document}

\title{Mirror Symmetry, Hitchin's Equations, And Langlands Duality}

\author{Edward Witten}
\address{School of Natural Sciences, Institute for Advanced Study, Princeton NJ 08540}
\email{dgaiotto@gmail.com, witten@ias.edu}
\thanks{Supported in part by NSF Grant Phy-0503584.}


\date{February, 2007}

\def\Bbb{\mathbb}

\begin{abstract}
Geometric Langlands duality can be understood from statements of
mirror symmetry that can be formulated in purely topological terms
for an oriented two-manifold $C$.  But understanding these
statements is extremely difficult without picking a complex
structure on $C$ and using Hitchin's equations.  We sketch the
essential statements both for the ``unramified'' case that $C$ is
a compact oriented two-manifold without boundary, and the
``ramified'' case that one allows punctures.  We also give a few
indications of why a more precise description requires a starting
point in four-dimensional gauge theory.
\end{abstract}

\maketitle

\input epsf
\section{The $A$-Model And The $B$-Model}
\label{intro}
\def\neg{\negthinspace}

Let $G$ be a compact Lie group and let  $G_\C$ be its
complexification. And let $C$ be  a compact oriented two-manifold
without boundary. We write $\cal Y(G,C)$ (or simply $\cal Y(G)$ or
$\cal Y$ if the context is clear) for the moduli space of flat
$G_\C$ bundles $E\to C$, modulo gauge transformations. Equivalently,
$\cal Y(G,C)$ parametrizes homomorphisms\footnote{Actually, it is
best to define $\cal Y(G,C)$ as a geometric invariant theory
quotient that parametrizes stable homomorphisms plus equivalence
classes of semi-stable ones.  This refinement will not concern us
here.  See  section \ref{belstacks}.} of the fundamental group of
$C$ to $G_\C$.

$\cal Y(G,C)$ is in a natural way a complex symplectic manifold,
that is a complex manifold with a nondegenerate holomorphic
two-form.  The complex structure comes simply from the complex
structure of $G_\C$, and the symplectic form, which we call
$\Omega$, comes from the intersection pairing\footnote{The
definition of this intersection pairing depends on the choice of
an invariant quadratic form on the Lie algebra of $G$.  It can be
shown using Hitchin's $\C^*$ action on the moduli space of Higgs
bundles that the $A$-model that we define shortly is independent
of this choice, up to a natural isomorphism.  The geometric
Langlands duality that one ultimately defines likewise does not
depend on this choice.} on $H^1(C,{\rm ad}(E))$, where
${\rm{ad}}(E)$ is the adjoint bundle associated to a flat bundle
$E$. Since $\cal Y(G,C)$ is a complex symplectic manifold, in
particular it follows that its canonical line bundle is naturally
trivial.

Geometric Langlands duality is concerned with certain topological
field theories associated with $\cal Y(G,C)$. The most basic of
these are the $B$-model that is defined by viewing $\cal Y(G,C)$ as
a complex manifold with trivial canonical bundle, and the $A$-model
that is defined by viewing it as a real symplectic manifold with
symplectic form\footnote{The usual definition of $\Omega$ is such
that ${\rm Im}\,\Omega$ is cohomologically trivial, while ${\rm
Re}\,\Omega$ is not.  The fact that $\omega={\rm Im}\,\Omega$ is
cohomologically trivial is a partial explanation of the fact,
mentioned in the last footnote, that the $A$-model of $\cal Y$ is
invariant under scaling of $\omega$. } $\omega={\rm Im}\,\Omega$.

These are the topological field theories that are relevant to the
most basic form of geometric Langlands duality.  However, there is
also a generalization that is relevant to what is sometimes called
quantum geometric Langlands.  {}From the $A$-model side, it is
obvious that a generalization is possible, since we could use a more
general linear combination of ${\rm Re}\,\Omega$ and ${\rm
Im}\,\Omega$ in defining the $A$-model.  What is less evident is
that the $B$-model can actually be deformed, as a topological field
theory, into this family of $A$-models.  This rather surprising fact
is natural from the point of view of generalized complex geometry,
see \cite{H1}, and has been explained from that point of view in
 section 4.6 of \cite{G}, as a general statement about complex
symplectic manifolds.  In \cite{KW}, sections 5.2 and 11.3, it was
shown that quantum geometric Langlands is naturally understood in
precisely this setting.

\def\Z{\Bbb{Z}}
Here, however, to keep things simple, we will focus on the most
basic $B$-model and $A$-model that were just described.

\section{Mirror Symmetry And Hitchin's Equations}

The next ingredient we need is Langlands or Goddard-Nuyts-Olive
duality. To every compact Lie group $G$ is naturally associated its
dual group $^L\neg G$.  For example, the dual of $SU(N)$ is
$PSU(N)=SU(N)/\Z_N$, the dual of $E_8$ is $E_8$, and so on.  And we
must also recall the concept of mirror symmetry between $A$-models
and $B$-models (for example, see \cite{KH}).  This is a quantum
symmetry of two-dimensional nonlinear sigma models whose most basic
role is to transform questions of complex geometry into questions of
symplectic geometry.  The geometric Langlands correspondence does
not appear at first sight to be an example of mirror symmetry, but
it turns out that it is.

With a little bit of hindsight (the question was first addressed
in \cite{HT}, following earlier work by \cite{BJV}, and
\cite{HMS}), we may ask whether the $B$-model of $\cal Y({}^L\neg
G,C)$ may be mirror to the $A$-model of $\cal Y(G,C)$. Even once
this question is asked, it is difficult to answer it without some
additional structure.  The additional structure that comes in
handy is provided by Hitchin's equations, see \cite{H2}.  Until
this point, $C$ has simply been an oriented two-manifold (compact
and without boundary). But now we pick a complex structure and
view $C$ as a complex Riemann surface. Hitchin's equations with
gauge group $G$ are equations for a pair $(A,\phi)$. Here $A$ is a
connection  on a $G$-bundle $E\to C$ (we stress that the structure
group of $E$ is now the {\it compact} group $G$), and
 $\phi$  is a one-form on $C$ with values in ${\rm ad}(E)$.
Hitchin's equations, which are elliptic modulo the gauge group,
are the system:
\begin{align}\label{hitchin}
F-\phi\wedge\phi &= 0 \nonumber \\
D\phi = D\star \phi &= 0. \\ \nonumber \end{align} Here $\star$ is
the Hodge star operator determined by the complex structure on $C$.
The role of the complex structure of $C$ is that it enables us to
write the last of these equations.

\def\A{{\cal A}}
\def\F{{\cal F}}
\def\Y{{\cal Y}}
A solution of Hitchin's equations has two interpretations.  On the
one hand, given such a solution, we can define the complex-valued
connection $\A=A+i\phi$.  Hitchin's equations imply that the
corresponding curvature $\F=d\A+\A\wedge \A$ vanishes, so a solution
of Hitchin's equations defines a complex-valued flat connection, and
thus a point in $\Y(G,C)$.

On the other hand, the $(0,1)$ part of the connection $A$ determines
a $\bar\partial$ operator on the bundle $E$ (or rather its
complexification, which we also call $E$).  There is no
integrability condition on $\bar\partial$ operators in complex
dimension 1, so this $\bar\partial $ operator endows $E$ with a
complex structure; it becomes a holomorphic $G_\C$ bundle over $C$.
Moreover, let us write $\phi=\varphi+\bar\varphi$, where $\varphi$
and $\bar\varphi$ are the $(1,0)$ and $(0,1)$ parts of $\phi$. Then
Hitchin's equations imply that $\varphi$, regarded as a section of
$K\otimes {\rm ad}(E)$ (with $K$ the canonical line bundle of $C$),
is holomorphic.  The pair $(E,\varphi)$, where $E\to C$ is a
holomorphic $G_\C$ bundle and $\varphi\in H^0(C,K\otimes{\rm
ad}(E))$, is known as a Higgs bundle.

\def\MH{{\cal M}_H}
We write $\MH$ for the moduli space of solutions of Hitchin's
equations, modulo a gauge transformation.  The fact that a solution
of these equations can be interpreted in two different ways means
that $\MH$ is endowed with two different natural complex structures.
In one complex structure, which has been called $I$, $\MH$
parametrizes isomorphism classes of semistable Higgs bundles
$(E,\varphi)$. In another complex structure, $J$, it parametrizes
equivalence classes of flat $G_\C$-bundles or in other words
homomorphisms $\rho:\pi_1(C)\to G_\C$.   $I$, $J$, and $K=IJ$ fit
together to a natural hyper-Kahler structure on $\MH$, as described
in \cite{H2}.  In particular, there are holomophic two-forms
$\Omega_I,\Omega_J,\Omega_K$ and Kahler forms
$\omega_I,\omega_J,\omega_K$.  These are all related by
$\Omega_I=\omega_J+i\omega_K$, and cyclic permutations of this
statement, as is usual in hyper-Kahler geometry.

In complex structure $J$, $\MH$ is the same as the variety $\cal Y$
that we described earlier.  The natural holomorphic symplectic form
$\Omega$ of $\cal Y$ is the same as $i\Omega_J$.  And the real
symplectic form $\omega={\rm Im}\,\Omega$ used in defining the
$A$-model coincides with $\omega_K$.  Complex structure $J$ and the
holomorphic symplectic form $\Omega_J=\omega_K+i\omega_I$ do not
depend on the chosen complex structure on $C$, in contrast to the
rest of the hyper-Kahler structure of $\MH$.

\remark As an aside, one may ask how closely related $\phi$, known
in the present context  as the Higgs field, is to the Higgs fields
of particle physics.  Thus, to what extent is the terminology that
was introduced in \cite{H2} actually justified?  The main
difference is that Higgs fields in particle physics are scalar
fields, while $\phi$ is a one-form on $C$ (valued in each case in
some representation of the gauge group). However, although
Hitchin's equations were first written down and studied directly,
they can be obtained from ${\cal N}=4$ supersymmetric gauge theory
via a sort of twisting procedure (similar to the procedure that
leads from ${\cal N}=2$ supersymmetric gauge theory to Donaldson
theory).  In this twisting procedure, some of the Higgs-like
scalar fields of ${\cal N}=4$ super Yang-Mills theory are indeed
converted into the Higgs field that enters in Hitchin's equations.
This gives a reasonable justification for the terminology.

\smallskip
As we will explain next, it is possible, with the aid of Hitchin's
equations,  to answer the question of whether the $B$-model of $\cal
Y({}^L\neg G,C)$ is mirror to the $A$-model of $\cal Y(G,C)$. This
in fact was first pointed out in \cite{HT}, and used in \cite{KW} as
a key ingredient in understanding the geometric Langlands
correspondence.

\section{The Hitchin Fibration}\label{hitchfib}

We will have to use the Hitchin fibration.  This is the map,
holomorphic in complex structure $I$, that takes a Higgs bundle
$(E,\varphi)$ to the characteristic polynomial of $\varphi$. For
example, for $G=SU(2)$, $(E,\varphi)$ is mapped simply to the
quadratic differential $\det\,\varphi$.  The target of the Hitchin
fibration is thus in this case the space $\cmmib B=H^0(C,K^2)$
that parametrizes quadratic differentials. This has a natural
analog for any $G$.  {}From the standpoint of complex structure
$I$, the generic fiber of the map $\pi:\MH\to \cmmib B$ is a
complex abelian variety (or to be slightly more precise, in
general a torsor for one). The fibers are Lagrangian from the
standpoint of the holomorphic symplectic form $\Omega_I$. Such a
fibration by complex Lagrangian tori turns $\MH$ into a completely
integrable Hamiltonian system in the complex sense  \cite{H3}.

There is, however, another way to look at the Hitchin fibration, as
first described in \cite{HT}. Let us go back to the $A$-model
defined with the real symplectic structure $\omega$. Since the
fibers of $\pi:\MH\to\cmmib B$ are Lagrangian for
$\Omega_I=\omega_J+i\omega_K$, they are in particular Lagrangian for
$\omega=\omega_K$.  Moreover, being holomorphic in complex structure
$I$, these fibers are actually area-minimizing in their homology
class -- here areas are computed using the hyper-Kahler metric on
$\MH$.  So the Hitchin fibration, from the standpoint of the
$A$-model,  is actually a fibration of $\MH$ by special Lagrangian
tori.

Mirror symmetry is believed to arise from $T$-duality on the fibers
of a special Lagrangian fibration, see \cite{SYZ}.  Generally, it is
very difficult to explicitly exhibit a non-trivial special
Lagrangian fibration.  The present example is one of the few
instances in which this can actually be done, with the aid of the
hyper-Kahler structure of $\MH$ and its integrable nature.
Non-trivial special Lagrangian fibrations are hard to understand
because it is difficult to elucidate the structure of the
singularities.  In the hyper-Kahler context, the fact that the
fibers are holomorphic in a different complex structure makes
everything far more accessible.

Once we actually find a special Lagrangian fibration, what we are
supposed to do with it, in order to give an example of mirror
symmetry, is to construct the dual special Lagrangian fibration,
which should be mirror to the original one. The mirror map
exchanges the symplectic structure on one side with the complex
structure on the other side.

In the present context, there is a very beautiful description of
the dual fibration: it is, as first shown in \cite{HT},  simply
the Hitchin fibration of the dual group.  Thus one considers
$\MH({}^L\neg G,C)$, the moduli space of solutions of Hitchin's
equation for the dual group $^L\neg G$. It turns out that the
bases of the Hitchin fibrations for $G$ and $^L\neg G$ can be
identified in a natural way.  The resulting picture is something
like this:
$$
\begin{array}{ccccc}
\MH({}^L\neg G,C) & \; & \; & \; & \MH(G,C) \\
\; & \searrow & \; & \swarrow & \; \\
\; & \; & {\cmmib B} & \; & \;
\end{array}
$$
In complex structure $I$, the fibers over a generic point $b\in
\cmmib B$ are, roughly speaking, dual abelian varieties (more
precisely, they are torsors for dual abelian varieties).

Alternatively, the fibers are special Lagrangian submanifolds in
the symplectic structure $\omega=\omega_K$.  {}From this second
point of view, the same picture leads to a mirror symmetry between
the $B$-model of $\MH({}^L\neg G,C)$ in complex structure $J$ and
the $A$-model of $\MH(G,C)$ with symplectic structure $\omega_K$.

As we have just explained, the tools that make this mirror
symmetry visible are the hyper-Kahler structure of $\MH$ and its
Hitchin fibration. Those structures depend on the choice of a
complex structure on $C$, but in fact, the resulting mirror
symmetry does not really depend on this choice. This was shown in
\cite{KW} in the process of deriving this example of mirror
symmetry from a four-dimensional topological field theory.  The
topological field theory in question is obtained by twisting of
${\cal N}=4$ super Yang-Mills theory.

\subsection{A Few Hints}

There are a few obstacles to overcome to go from this instance of
mirror symmetry to the usual formulation of geometric Langlands
duality.  Unfortunately, it will not be practical here to give more
than a few hints.

One ~ key ~~ point ~ is ~ that ~ in ~ the ~ usual ~ formulation, ~
the ~ dual ~ of ~ a ~ $B$-brane ~ on ~ \linebreak $\MH({}^L\neg
G,C)$ is supposed to be not an $A$-brane on $\MH(G,C)$ -- which is
what we most naturally get from the above construction -- but a
sheaf of ${\cal D}$-modules on $\M(G,C)$, the moduli space of
$G$-bundles over $C$ (a sheaf of ${\cal D}$-modules is by
definition a sheaf of modules for the sheaf ${\cal D}$ of
differential operators on $\M(G,C)$). The link between the two
statements is explained in \cite{KW}, section 11, using the
existence of a special $A$-brane on $\MH(G,C)$ that is intimately
related to differential operators on $\M(G,C)$. This relation is
possible because, as explained in \cite{H2}, $\MH(G,C)$ contains
$T^*\M^{st}(G,C)$ as a Zariski open set; here $\M^{st}(G,C)$ is
the subspace of $\M(G,C)$ parametrizing strictly stable bundles.

Another key point is the following. A central role in the usual
formulation is played by the geometric Hecke operators, which act on
holomorphic $G$-bundles over $C$ and therefore also on ${\cal
D}$-modules on $\M(G,C)$.  They have a natural role  in the present
story, but this is one place that one misses something if one
attempts to express this subject just in terms of two-dimensional
sigma models and mirror symmetry.  This particular instance of
mirror symmetry actually originates from a duality in an underlying
four-dimensional gauge theory. Once this is understood, basic facts
about the Wilson and 't Hooft line operators of gauge theory  lead
to the usual statements about Hecke eigensheaves, as explained in
sections 9 and 10 of \cite{KW}. The geometric Hecke operators are
naturally reinterpreted in this context in terms of the Bogomolny
equations of three-dimensional gauge theory, which are of great
geometrical as well as physical interest and have been much studied,
for example in \cite{AH}.

A proper formulation of some of these statements leads to another
important role for four dimensions.  The usual formulation of
geometric Langlands involves ${\cal D}$-modules not on the moduli
space of semistable $G$-bundles but on the moduli stack of all
$G$-bundles. The main reason for this is that one cannot see the
action of the Hecke operators if one considers only semistable
bundles. As we will explain in section \ref{stacks}, the role of
stacks in the standard description can be understood as a strong
clue for an alternative approach that starts in four-dimensional
gauge theory.

\section{Ramification}

\def\CC{{\eusm C}}
Before getting back to stacks, however, I want to give an idea of
what is called ``ramification'' in the context of geometric
Langlands.

A simple generalization of what we have said so far is to consider
flat bundles not on a closed oriented two-manifold $C$ but on a
punctured two-manifold $C'=C\backslash p$; that is, $C'$ is $C$ with
a point $p$ omitted.

We pick a conjugacy class $\CC\subset G_\C$, and we let
$\Y(G,C';\CC)$ denote the moduli space of homomorphisms
$\rho:\pi_1(C')\to G_\C$, up to conjugation, such that the monodromy
around $p$ is in the conjugacy class $\CC$.

Many statements that we made before have natural analogs in this
punctured case.  In particular, $\Y(G,C';\CC)$ has a natural
structure of a complex symplectic manifold.  It has a natural
complex structure  and holomorphic symplectic form $\Omega$. Just as
in the unpunctured case, we can define a $B$-model of
$\Y(G,C';\CC)$. Also, viewing $\Y(G,C';\CC)$ as a real symplectic
manifold with symplectic form $\omega={\rm Im}\,\Omega$, we can
define an $A$-model.  The $B$-model and the $A$-model are both
completely independent of the complex structure of $C'$.

Next, introduce the dual group $^L\neg G$ and let $^L\neg\CC$ denote
a conjugacy class in its complexification.  Again, the space
$\Y({}^L\neg G,C';{}^L\neg \CC)$ has a natural $B$-model and
$A$-model.

Based on what we have said so far, one may wonder if, for some map
between $\CC$ and $^L\neg \CC$, there might be a mirror symmetry
between $\Y(G,C';\CC)$ and $\Y({}^L\neg G,C';{}^L\neg \CC)$.  The
answer to this question is ``not quite,'' for a number of reasons.
One problem is that there is no natural correspondence between
conjugacy classes in $G_\C$ and in ${}^L\neg G_\C$.  A more
fundamental problem is that the $B$-model of $\Y(G,C';\CC)$ varies
holomorphically with the conjugacy class $\CC$, but the $A$-model of
the same space does not.  To find a version of the statement that
has a chance of being right, we have to add additional parameters to
find a mirror-symmetric set.

In any event, regardless of what parameters one adds, it is very
difficult to answer the question about mirror symmetry if $C'$ is
viewed simply as an oriented two-manifold with a puncture.  To make
progress, just as in the unramified case (that is, the case without
punctures), it is very helpful to endow $C'$ with a complex
structure and to use Hitchin's equations.  This actually also helps
us in finding the right parameters, because an improved set of
parameters appears just in trying to give a natural formulation of
Hitchin's equations on a punctured surface.  Let $z$ be a local
parameter near the puncture and write $z=re^{i\theta}$.  In the
punctured case, it is natural, see \cite{S},  to introduce variables
$\alpha,\beta,\gamma$ taking values in the Lie algebra $\mathfrak t$
of a maximal torus $T\subset G$, and consider solutions of Hitchin's
equations on $C'$ whose behavior near $z=0$ is as follows:
\begin{align}\label{curtsy}A &
=  \alpha \,d\theta +\dots\\
           \phi & = \beta\,{\frac{dr}{r}}-\gamma\,d\theta+\dots.\end{align}
The ellipses refer to terms that are less singular near $z=0$.

All the usual statements about Hitchin's equations have close
analogs in this situation. The moduli space of solutions of
Hitchin's equations with this sort of singularity is a hyper-Kahler
manifold $\MH(G,C;\alpha,\beta,\gamma)$.  In one complex structure,
usually called $J$, it coincides with $\Y(G,C;\CC)$, where $\CC$ is
the conjugacy class that contains\footnote{For simplicity, we assume
that $U$ is regular.  The more involved statement that holds in
general is explained in \cite{GW}.} $U=\exp(-2\pi(\alpha-i\gamma))$.
In another complex structure, often called $I$,
$\MH(G,C;\alpha,\beta,\gamma)$ parametrizes Higgs bundles
$(E,\varphi)$, where $\varphi\in H^0(C',K\otimes {\rm ad}(E))$ has a
pole at $z=0$, with $\varphi\sim \frac{1}{2}(\beta+i\gamma)(dz/z)$.
Moreover, there is a Hitchin fibration, and most of the usual
statements about the unramified case -- those that we have explained
and those that we have omitted here -- have close analogs.  For a
much more detailed explanation, and references to the original
literature, see \cite{GW}.

The variables $\alpha,\beta,\gamma$ are a natural set of
parameters for the classical geometry.  However, quantum
mechanically, there is one more natural variable $\eta$ (analogous
to the usual $\theta$-angles of gauge theory), as described in
section 2.3 of \cite{GW}. With the complete set of parameters
$(\alpha,\beta,\gamma,\eta)$ at hand, it is possible to formulate
a natural duality statement, according to which $\MH({}^L\neg
G,C;{}^L\neg\alpha,{}^L\neg\beta,{}^L\neg \gamma,{}^L\neg \eta)$
is mirror to $\MH(G,C;\alpha,\beta,\gamma,\eta)$, with a certain
map between the parameters, described in section 2.4 of \cite{GW}.
The main point of this map is that
$(\alpha,\eta)=({}^L\eta,-{}^L\alpha)$.  Since the monodromy $U$
depends on ${}^L\neg\alpha$, this shows that the dual of the
monodromy involves the quantum parameter $\eta$ that is invisible
in the classical geometry.  In the $A$-model, $\eta$ becomes the
imaginary part of the complexified Kahler class.

This duality statement leads, after again mapping $A$-branes to
${\cal D}$-modules, to a statement of geometric Langlands duality
for this situation similar to what has been obtained via algebraic
geometry and two-dimensional conformal field theory in \cite{FG}.

\remark We pause here to explain one very elementary fact about the
classical geometry that will be helpful as background for section
\ref{wild}.  In complex structure $J$, a solution of Hitchin's
equations with the singularity of eqn. (\ref{curtsy}) describes a
flat $G_\C$ bundle $E\to C'$ with monodromy around the puncture $p$.
$E$ can be extended over $p$ as a holomorphic bundle, though of
course not as a flat one, and moreover from a holomorphic point of
view, $E$ can be trivialized near $p$.  The flat connection on $E\to
C'$ is then represented, in this gauge, by a holomorphic
$(1,0)$-form on $C'$ (valued in the Lie algebra of $G_\C$) with a
simple pole at $p$:
\begin{equation}\label{zork}\cal
A=dz\left(\frac{\alpha-i\gamma}{iz}+\dots\right),\end{equation}
where the omitted terms are regular at $z=0$.  The singularity of
the connection at $z=0$ is a simple pole because the ansatz
(\ref{curtsy}) for Hitchin's equations only allows a singularity of
order $1/|z|$.  A holomorphic connection with such a simple pole is
said to have a regular singularity.

In geometric Langlands, what is usually called tame ramification is,
roughly speaking, the case that we have just arrived at: a
holomorphic bundle $E\to C$ that has a holomorphic connection form
with a regular singularity. Actually, the phrase ``tame
ramification'' is sometimes taken to refer to the case that the
residue of the simple pole is nilpotent, while in eqn. (\ref{zork})
we seem to be in the opposite case of semi-simple residue. In
\cite{GW}, it is explained that, with some care, mirror symmetry for
$\MH(G,C';\alpha,\beta,\gamma,\eta)$ is actually a sufficient
framework to understand geometric Langlands for a connection with a
simple pole of any residue.  For example, the case of a nilpotent
residue can be understood by setting
${}^L\neg\alpha={}^L\neg\gamma=0$ (or $\gamma=\eta=0$ in the dual
description).

\section{Wild Ramification}\label{wild}

Based on an analogy with number theory, geometric Langlands is
usually formulated not only for the case of tame ramification.  One
goes on to inquire about an analogous duality statement involving a
holomorphic bundle $E\to C$ with a holomorphic connection that has a
pole of any order.  In other words, after trivializing the
holomorphic structure of $E$ near a point $p\in C$, the connection
looks like
\begin{equation}\label{gelg}\cal
A=dz\left(\frac{T_n}{z^n}+\frac{T_{n-1}}{z^{n-1}}+\dots+\frac{T_1}{z}+\dots\right),\end{equation}
where regular terms are omitted.  A meromorphic connection with a
pole of degree greater than 1 is said to have an irregular
singularity.

Trying to formulate a duality statement for this situation poses, at
first sight, a severe challenge for the approach to geometric
Langlands  described here.  Our basic point of view is that the
fundamental duality statements depend on $C$ only as an oriented
two-manifold. A complex structure on $C$ is introduced only as a
tool to answer certain natural questions that can be asked without
introducing the complex structure.

{}From this point of view, tame ramification is natural because a
simple pole in this sense has a clear topological meaning.  A
meromorphic connection with a simple pole at a point $p\in C$ is a
natural way to encode the monodromy about $p$ of a flat connection
on $C'=C\backslash p$.  And this monodromy, of course, is a purely
topological notion.   But what could possibly be the topological
meaning of a connection with a pole of degree greater than 1?

A closely related observation is that $T_1$ is the residue of the
pole in $\cal A$ at $z=0$, and so is independent of the choice of
local coordinate $z$. However, the coefficients $T_2,\dots,T_n$ of
the higher order poles most definitely do depend on the choice of a
local coordinate. How can we hope to include them in a theory that
is supposed to depend on $C$ only as an oriented two-manifold?

Moreover, if the plan is to formulate a duality conjecture of a
topological nature and then prove it using Hitchin's equations, we
face the question of whether Hitchin's equations are compatible with
an irregular singularity.  Hitchin's equations for a pair
$\Phi=(A,\phi)$ are schematically of the form $d\Phi+\Phi^2=0$.  If
near $z=0$, we have a singularity with $|\Phi|\sim 1/|z|^n$, then
$|d\Phi|\sim 1/|z|^{n+1}$ and $|\Phi|^2\sim 1/|z|^{2n}$.  For $n=1$,
$d\Phi$ and $|\Phi|^2$ are comparable in magnitude, and therefore
Hitchin's equations look reasonable. However, for $n>1$, we have
$|\Phi|^2>>|d\Phi|$, and it looks like the nonlinear term in
Hitchin's equations will be too strong.

Both questions, however, have natural answers.  The answer to the
first question is that, despite appearances, one actually can
associate to a connection with irregular singularity something that
goes beyond the ordinary monodromy and has a purely topological
meaning. The appropriate concept is an extended monodromy that
includes Stokes matrices as well as the ordinary monodromy.  Stokes
matrices are part of the classical theory of ordinary differential
equations with irregular singularity; for example, see \cite{Wa}.

Assuming for brevity that the leading coefficient $T_n$ of the
singular part of the connection is regular semi-simple, one can
make a gauge transformation to conjugate $T_1,\dots,T_n$ to the
maximal torus. Then one defines a moduli space ${\cal
Y}(G,C;T_1,\dots,T_n)$ that parametrizes, up to a gauge
transformation, pairs consisting of a holomorphic $G_\C$-bundle
over $C$ and a connection with an irregular singularity of the
form described in eqn. (\ref{gelg}). As shown in \cite{B}, it
turns out that this space ${\cal Y}(G,C;T_1,\dots,T_n)$ is in a
natural way a complex symplectic manifold, with a complex
symplectic structure that depends on $C$ only as an oriented
two-manifold. This can be proved by adapting to the present
setting the gauge theory definition of the symplectic structure,
formulated in \cite{ABott}. Moreover the complex symplectic
structure of ${\cal Y}(G,C;T_1,\dots,T_n)$ is independent of
$T_2,\dots,T_n$ (as long as $T_n$ remains semi-simple).

At  this  point  the  important concept of isomonodromic
deformation, introduced by
 \cite{JMU}, comes into play.  There  is  a natural
way to vary the parameters  $T_2,\dots,T_n$, without ~~~ changing
~~~ the ~~~ generalized  ~~~ monodromy ~~~ data ~~~ that ~~~ is
parametrized ~~~ by ~~~  ${\cal Y}(G,C;T_1,\dots,T_n)$. Moreover,
as has been proved quite recently in
 \cite{B},  the  complex symplectic   structure   of the
space of  generalized  monodromy  data  is  invariant  under
isomonodromic deformation.  Thus, roughly speaking,  one  can define
a complex  symplectic  manifold
 ${\cal Y}(G,C;T_1,n)$ that
depends only on $T_1$ and the integer $n\geq 1$.

The fact that the parameters $T_2,\dots,T_n$ turn out to be, in the
sense just described, inessential, is certainly welcome, since as we
have already observed, they have no evident topological meaning. Now
we are in a situation very similar to what we had in the unramified
and tamely ramified cases.  Given ${\cal Y}(G,C;T_1,n)$ as a complex
symplectic manifold, with complex symplectic form $\Omega$, we can
define its $B$-model, or its $A$-model using the real symplectic
form $\omega={\rm Im}\,\Omega$.  Of course, we can do the same for
the dual group, defining another complex symplectic manifold {${\cal
Y}({}^L\neg G,C;{}^LT_1,n)$}, with its own $B$-model and $A$-model.
And, just as in the unramified case, we can ask if these two models
are mirror to each other.

Even before trying to answer this question, we should refine it
slightly.  Because of the constraint that $T_n$ should be regular
semi-simple, it is not quite correct to simply forget about $T_n$.
There can be monodromies as $T_n$ varies.  We think of $T_n$ as
taking values in $\mathfrak t_\C^{\rm reg}\otimes K_p^{n-1}$, with
notation as follows: $\mathfrak t_\C$ is the Lie algebra of a
maximal torus in $G_\C$, $\mathfrak t_\C^{\rm reg}$ is its subspace
consisting of regular elements, and $K_p$ is the fiber at $p$ of the
cotangent bundle to $C$. The fundamental group of $\mathfrak
t_\C^{\rm reg}$ is known as the braid group of $G$; we call it
$B(G)$. Because of the monodromies, one really needs to choose a
basepoint $*\in \mathfrak t_\C^{\rm reg}$ to define ${\cal
Y}({}^L\neg G,C;{}^LT_1,n)$; to be more precise, we can denote this
space as ${\cal Y}({}^L\neg G,C;{}^LT_1,n,*)$. The group $B(G)$ acts
via monodromies on both the $B$-model and the $A$-model of ${\cal
Y}(G,C;T_1,n,*)$. Dually, the corresponding braid group $B({}^L\neg
G)$ acts on the $B$-model and the $A$-model of ${\cal Y}({}^L\neg
G,C;{}^LT_1,n,*)$. However, the two groups $B(G)$ and $B({}^L\neg
G)$ are naturally isomorphic; indeed, modulo a choice of an
invariant quadratic form, there is a natural map from $\mathfrak
t_\C^{\rm reg}$ to $^L\neg \mathfrak t_\C^{\rm reg}$, so the two
spaces have the same fundamental group and a choice of basepoint in
one determines a basepoint in the other, up to homotopy. A better
(but still not yet precise) question is whether there is a mirror
symmetry between ${\cal Y}(G,C;T_1,n,*)$ and ${\cal Y}({}^L\neg
G,C;{}^LT_1,n,*)$ that commutes with the braid group.

We expect as well that this mirror symmetry depends on $C$ only as
an oriented two-manifold, and so commutes with the mapping class
group. We can think of the mapping class group of $C$ and the braid
group as playing quite parallel roles.  In fact, because of the
appearance of the fiber $K_p$ of the canonical bundle in the last
paragraph, these two groups do not simply commute with each other;
the group that acts is an extension of the mapping class group by
$B(G)$.

Just as in the tamely ramified case, to get the right mirror
symmetry conjecture, we need to extend the parameters slightly to
get a mirror-symmetric set.  But we also face the fundamental
question of whether Hitchin's equations are compatible with wild
ramification. As explained above, the nonlinearity of Hitchin's
equations makes this appear doubtful at first sight. But happily,
it turns out that all is well, as shown in \cite{BB}. The key
point is that, again with $T_n$ assumed to be regular semi-simple,
we can assume that the singular part of the connection is abelian.
Though Hitchin's equations are nonlinear, they become linear in
the abelian case, and once abelianized, they are compatible with a
singularity of any order. Using this as a starting point, it turns
out that, for any $n$, one can develop a theory of Hitchin's
equations with irregular singularity that is quite parallel to the
more familiar story in the unramified case. For example, the
moduli space $\MH$ of solutions of the equations is hyper-Kahler.
In one complex structure, $\MH$ parametrizes flat connections with
a singularity similar to that in eqn. (\ref{gelg}); in another
complex structure, it parametrizes Higgs bundles $(E,\varphi)$ in
which $\varphi$ has an analogous pole of order $n$. There is a
Hitchin fibration, and all the usual properties have close
analogs.

All this gives precisely the right ingredients to use Hitchin's
equations to establish the desired mirror symmetry between the two
moduli spaces. See \cite{W} for a detailed explanation in which this
classical geometry is embedded in four-dimensional gauge theory.
 Many of the arguments are quite similar to those given in the tame case
in \cite{GW}. The construction makes it apparent that the duality
commutes with isomonodromic deformation.

Finally, one might worry that the assumption that $T_n$ is regular
semi-simple may have simplified things in some unrealistic way.
This is actually not the case.  For one thing, the analysis in
\cite{BB} requires only that $T_2,\dots,T_n$ should be
simultaneously diagonalizable (in some gauge), and in particular
semi-simple, but not that $T_n$ is regular. But even if these
coefficients are not semi-simple, there is no essential problem.
In the classical theory of ordinary differential equations, it is
shown that given any such equation with an irregular singularity
at $z=0$, after possibly passing to a finite cover of the
punctured $z$-plane and changing the extension of a holomorphic
bundle over the puncture at $z=0$, one can reduce to the case that
the irregular part of the singularity has the properties assumed
in \cite{BB}. Given this, one can adapt all the relevant arguments
concerning geometric Langlands duality to the more general case,
as is explained in section 6 of \cite{W}.

\section{Four-Dimensional Gauge Theory And Stacks}\label{stacks}

To a physicist, it is natural, in studying dualities involving gauge
theory, to begin in four dimensions, which is often found to be the
natural setting for gauge theory duality. There is a simple reason
for this.  The curvature, which is one of the most fundamental
concepts in gauge theory, is a 2-form.  In $d$ dimensions, the dual
of a $2$-form is a $(d-2)$-form, so it is only a 2-form if $d=4$.
This suggests that $d=4$ is the most natural dimension in which the
dual of a gauge theory might be another gauge theory.

Moreover, $\cal N=4$ supersymmetric Yang-Mills theory, originally
constructed in \cite{BSS}, is a natural place to start, as it has
the maximal possible supersymmetry, and has the celebrated duality
whose origins go back to the  early work of  \cite{MO}.  It indeed
turns out that geometric Langlands has a natural origin in a
twisted version of $\cal N=4$ super Yang-Mills theory in four
dimensions. The twisting is quite analogous to the twisting of
$\cal N=2$ super Yang-Mills theory that leads to Donaldson theory.

That particular motivation may seem opaque to some, and instead I
will adopt here a different approach in explaining why it is natural
to begin in four dimensions for understanding geometric Langlands,
instead of relying only on the $B$-model and $A$-model of
$\MH(G,C)$.

First of all, the $B$-model and the $A$-model of any space $X$ are
both twisted versions of a quantum sigma model that governs maps
$\Phi:\Sigma\to X$, where $\Sigma$ is a two-manifold (or better, a
supermanifold of bosonic dimension two). Since the $A$-model
involves in its most elementary form a counting of holomorphic
maps $\Phi:\Sigma\to X$ that obey appropriate conditions, the
roles of $\Sigma$ and $\Phi$ are clear in the $A$-model.  Mirror
symmetry indicates that it must be correct to also formulate the
$B$-model in terms of maps $\Phi:\Sigma\to X$, and this is done in
the usual formulation by physicists.

In the present case, we are interested, roughly speaking, in the
$B$- and $A$-models of $\MH(G,C)$, for some compact Lie group $G$
and two-manifold $C$.  Therefore, roughly speaking, we want to study
a sigma model of maps $\Phi:\Sigma\to \MH(G,C)$, where as before
$\Sigma$ is an auxiliary two-manifold.

The reason that this description is rough is that $\MH(G,C)$ has
singularities,\footnote{Moreover, these singularities are worse than
orbifold singularities.  Orbifold singularities would cause no
difficulty.  See \cite{FW} for a discussion of orbifold
singularities in geometric Langlands.} and the sigma model of target
$\MH(G,C)$ is therefore not really well-defined. Therefore a
complete description cannot be made purely in terms of a sigma model
in which the target space is $\MH(G,C)$, viewed as an abstract
manifold.  We need a more complete description that will tell us how
to treat the singularities.  What might this be?

By definition, a point in $\MH(G,C)$ determines up to
gauge-equivalence a pair $(A,\phi)$ obeying Hitchin's equations. $A$
and $\phi$ are fields defined on $C$, so let us write them more
explicitly as $(A(y),\phi(y))$, where $y$ is a coordinate on $C$.

Now suppose that we have a map $\Phi:\Sigma\to \MH(G,C)$, where
$\Sigma$ is a Riemann surface with a local coordinate $x$.  Such a
map is described by a pair  $(A(y),\phi(y))$ that also depends on
$x$. So we can describe the map $\Phi$ via fields
$(A(x,y),\phi(x,y))$ that depend on both $x$ and $y$.  We would like
to interpret these fields as fields on the four-manifold
$M=\Sigma\times C$.   The pair $(A(x,y),\phi(x,y))$ is not quite a
natural set of fields on $M$ but can be naturally completed to one.
For example, $A(x,y)$ is locally a one-form tangent to the second
factor in $M=\Sigma\times C$; to get a four-dimensional gauge field,
we should relax the condition that $A$ is tangent to the second
factor. Similarly, we can extend $\phi$ to an adjoint-valued
one-form on $\Sigma\times C$. ${\cal N}=4$ super Yang-Mills theory,
or rather its twisted version that is related to geometric
Langlands, is obtained by completing this set of fields to a
supersymmetric combination in a minimal fashion.

In ${\cal N}=4$ super Yang-Mills theory, there are no singularities
analogous to the singularities of $\MH(G,C)$.  The space of gauge
fields, for example, is an affine space, and  the other fields (such
as $\phi$) take values in linear spaces.  The problems with
singularities that make it difficult to define a sigma model of maps
$\Phi:\Sigma\to \MH(G,C)$ have no  analog in defining gauge theory
on $M=\Sigma\times C$ (or any other four-manifold).  The relation
between the two is that the two-dimensional sigma model is an
approximation to the four-dimensional gauge theory.  The
approximation breaks down when one runs into the singularities of
$\MH(G,C)$.  Any question that involves those singularities should
be addressed in the underlying four-dimensional gauge theory.

But away from singularities, it suffices to consider only the
smaller set of fields that describe a map $\Phi:\Sigma\to \MH(G,C)$.
Many questions do not depend on the singularities and for these
questions the description via two-dimensional sigma models and
mirror symmetry is adequate.

\subsection{Stacks}\label{belstacks}

To conclude, we will make contact with the counterpart of this
discussion in the usual mathematical theory.  We start with bundles
rather than Higgs bundles because this case will be easier to
explain.

In the usual mathematical theory, the right hand side of the
geometric Langlands correspondence is  described in terms of ${\cal
D}$-modules on, roughly speaking, the moduli space of all
holomorphic $G_\C$ bundles on the Riemann surface $C$.

However, instead of the moduli space $\M(G,C)$ of semi-stable
holomorphic $G_\C$ bundles  $E\to C$, one considers $\cal D$-modules
on the ``stack'' ${\rm Bun}_G(C)$ of all such bundles.  The main
reason for this is that to define the action of Hecke operators, it
is necessary to allow unstable bundles.  Unstable bundles are
related to the non-orbifold singularities of $\M(G,C)$.

What is a stack?  Roughly, it is a space that can everywhere be
locally described as a quotient.  The trivial case is a stack that
can actually be described globally as a quotient. Interpreting ${\rm
Bun}_G(C)$ as a global quotient would mean finding a pair
$(Y,W_\C)$, consisting of a smooth algebraic variety $Y$ and a
complex Lie group $W_\C$ acting on $Y$, with the following
properties.  Isomorphism classes of holomorphic $G_\C$ bundles $E\to
C$ should be in 1-1 correspondence with $W_\C$ orbits on $Y$, and
for every $E\to C$, its automorphism group should be isomorphic to
the subgroup of $W_\C$ leaving fixed the corresponding point in $Y$.

A pair $(Y,W_\C)$ representing in this way the stack ${\rm
Bun}_G(C)$ does not exist if $Y$ and $W_\C$ are supposed to be
finite-dimensional.  Indeed, the $G_\C$-bundle $E\to C$ can be
arbtirarily unstable, so there is no upper bound on the dimension
of its automorphism group.  So no finite-dimensional $W_\C$ can
contain all such automorphism groups as subgroups.

However as shown in \cite{ABott}, taking $G$ to be of adjoint type
for simplicity, there is a natural infinite-dimensional pair
$(Y,W_\C)$. One simply takes $Y$ to be the space of all
connections on a given $G$-bundle $E\to C$ which initially is
defined only topologically. One defines $W$ to be the group of all
gauge transformations of the bundle $E$; thus, if $E$ is
topologically trivial, we can identify $W$ as the group ${\rm
Maps}(C,G)$.  Then we take $W_\C$ to be the complexification of
$W$, or in other words ${\rm Maps}(C,G_\C)$. (This
complexification acts on $Y$ as follows. We associate to a
connection $A$ the corresponding $\bar\partial$ operator
$\bar\partial_A$. Then a complex-valued gauge transformation acts
by $\bar\partial_A\to g\bar\partial_A g^{-1}$.)

Suppose then that we were presented with the problem of making sense
of the supersymmetric sigma model of maps $\Phi:\Sigma\to \M(G,C)$,
given the singularities of $\M(G,C)$. (This is a practice case for
our actual problem, which involves $\MH(G,C)$ rather than
$\M(G,C)$.) Our friends in algebraic geometry would tell us to
replace $\M(G,C)$ by the stack ${\rm Bun}_G(C)$. We interpret this
stack as the pair $(Y,W_\C)$, where $Y$ is the space of all
connections on a $G$-bundle $E\to C$  and $W_\C$ is the complexified
group of gauge transformations. The connected components of the
stack correspond to the topological choices for $E$.

By a supersymmetric sigma model with target a pair $(Y,W_\C)$, with
$W_\C$ a complex Lie group acting on a complex manifold $Y$, we
mean\footnote{For a discussion of this construction in relation to
stacks, see \cite{PS}.} in the finite-dimensional case a
gauge-invariant supersymmetric sigma model in which the gauge group
is $W$ (a maximal compact subgroup of $W_\C$) and the target is $Y$.
Actually, to define this sigma model, we want $Y$ to be a Kahler
manifold with an $W$-invariant (but of course not $W_\C$-invariant)
Kahler structure.  The sigma model action contains a term which is
the square of the moment map for the action of $W$. This term is
minimized precisely when the moment map vanishes.  The combined
operation of setting the moment map to zero and dividing by $W$ is
equivalent classically to dividing by $W_\C$.

To write down the term in the action that involves the square of the
moment map (and in fact, to write down the kinetic energy of the
gauge fields) one needs an invariant and positive definite quadratic
form on the Lie algebra of $W$. If $W$ is finite-dimensional,
existence of such a form is equivalent to $W$ being compact.
However, the appropriate quadratic form also exists in the
infinite-dimensional case that $W={\rm Maps}(C,G)$ for some space
$C$. (An element of the Lie algebra of $W$ is a $\mathfrak g$-valued
function $\epsilon$ on $C$, and the quadratic form is defined by
$\int_C d\mu \,(\epsilon,\epsilon)$, where $(~,~)$ is an invariant
positive-definite quadratic form on $\mathfrak g$, and $\mu$ is a
suitable measure on $C$.)

Now, suppose we construct the two-dimensional sigma model of maps
from a Riemann surface $\Sigma$ to the stack ${\rm Bun}_G(C)$,
understood as above. What is the group of gauge transformations of
the sigma model? In general, in a gauge theory on any space
$\Sigma$ with gauge group $W$, the group of gauge transformations
(of a topologically trivial $W$-bundle, for simplicity) is the
group of maps from $\Sigma$ to $W$, or ${\rm Maps}(\Sigma,W)$. In
our case, $W$ is in turn ${\rm Maps}(C,G)$.  So ${\rm
Maps}(\Sigma,W)$ is the same as ${\rm Maps}(M,G)$, where
$M=\Sigma\times C$. But this is simply the group of gauge
transformations in gauge theory on $M$, with gauge group $G$. In
the present case, $\Sigma$ and $C$ are two-manifolds and $M$ is a
four-manifold. We have arrived at four-dimensional gauge theory.
If we chase through the definitions a little more, we learn that
the supersymmetric sigma model of maps $\Phi:\Sigma \to {\rm
Bun}_G(C)$ should be understood as four-dimensional ${\cal N}=2$
supersymmetric gauge theory, with gauge group $G$, on the
four-manifold $M=\Sigma\times C$.  (This is the theory that after
twisting is related to Donaldson theory.)

Now let us return to the original problem. Geometric Langlands
duality is a statement about the $B$-model and $A$-model not of
$\M(G,C)$ but of $\MH(G,C)$, the corresponding moduli space of
Higgs bundles, and its analog for the dual group $^L\neg G$.  To
deal with the singularities, we want to ``stackify'' this
situation.   We are now in  a hyper-Kahler context and the
appropriate concept of a stack should incorporate this.  (What
algebraic geometers would call the stack of Higgs bundles does not
quite do justice to the situation, since it emphasizes one complex
structure too much.) Since quaternionic Lie groups do not exist,
we cannot ask to construct $\MH(G,C)$ as the quotient of a smooth
space by a quaternionic Lie group.  However, the notion of a
symplectic quotient does have a good analog in the hyper-Kahler
world, namely the hyper-Kahler quotient, described in \cite{H4}.
 The
analog of what we explained for $\M(G,C)$ is to realize $\MH(G,C)$
as the hyper-Kahler quotient  of a smooth space $Y$ by a group
$W$.
 It may be impossible to do this with finite-dimensional $Y$ and
$W$, but in the infinite-dimensional world, this problem has a
natural solution described in Hitchin's original paper on the
Hitchin equations, \cite{H2}. ($Y$ is the space of pairs $(A,\phi)$
on $C$, and $W={\rm Maps}(C,G)$.) Taking this as input and
interpreting what it should mean to have a sigma model whose target
is the hyper-Kahler stack corresponding to $\MH(G,C)$, one arrives
at the twisted version of $\cal N=4$ super Yang-Mills theory that
was the starting point in \cite{KW}.

\bibliographystyle{agsm}
\bibliography{Hitchin}

\end{document}